\newcommand{\commentt}[2]{#1}
\newcommand{\comment}[1]{}
\newcommand{\titt}{Optimal unbiased estimators via convex hulls}

\newcommand{\var}{{\rm Var}}
\newcommand{\std}{{\rm Std}}
\commentt{
\newcommand{\citet}{\citeasnoun}
\documentclass[a4paper,11pt]{article}
\usepackage{graphicx,slashbox,amsthm}
\usepackage[full]{harvard}

\usepackage[ps,dvips]{xy}
\usepackage{pgfplots}
\title{\titt}
\author{Nabil Kahal\'e
\thanks{\emph{ESCP Europe and  Labex ReFi, 75011 Paris,
France; {e-mail: }{nkahale@escpeurope.eu}.}}
}
\date{\today}
\usepackage{amstext}
\usepackage{amssymb}
\usepackage{amsmath}
\usepackage{setspace}
\usepackage[margin=1in]{geometry}
\usepackage[english]{babel}
\bibliographystyle{dcu}
\usepackage{color}
\usepackage{textcomp}
\usepackage{varioref}

\begin{document}

\newtheorem{example}{Example}[section]
\newtheorem{theorem}{Theorem}[section]
\newtheorem{conjecture}{Conjecture}[section]
\newtheorem{lemma}{Lemma}[section]
\newtheorem{proposition}{Proposition}[section]
\newtheorem{remark}{Remark}[section]
\newtheorem{corollary}{Corollary}[section]
\newtheorem{definition}{Definition}[section]
\numberwithin{equation}{section}
\maketitle
\newcommand{\ABSTRACT}[1]{\begin{abstract}#1\end{abstract}}
\newcommand{\citep}{\cite}
}
{
\documentclass[opre,nonblindrev]{informs3} 
\usepackage[margin=1in]{geometry}
\usepackage{textcomp}
\usepackage{pgfplots}
\DoubleSpacedXII



\usepackage[round]{natbib}
 \bibpunct[, ]{(}{)}{,}{a}{}{,}%
 \def\bibfont{\small}%
 \def\bibsep{\smallskipamount}%
 \def\bibhang{24pt}%
 \def\newblock{\ }%
 \def\BIBand{and}%

\TheoremsNumberedThrough     
\ECRepeatTheorems

\EquationsNumberedThrough    

\MANUSCRIPTNO{} 

\renewcommand{\qed}{\halmos}
\renewcommand{\cite}{\citep}

\begin{document}
\RUNAUTHOR{Kahal\'e}
\RUNTITLE{Multilevel  methods for discrete Asian options}
\TITLE{\titt}
\ARTICLEAUTHORS{%
\AUTHOR{Nabil Kahal\'e}

\AFF{ESCP Europe and Labex Refi, 75011 Paris, France, \EMAIL{nkahale@escpeurope.eu} \URL{}}
} 
\bibliographystyle{informs2014} 
}
\ABSTRACT{
Necessary and sufficient conditions for the square-integrability of recently proposed unbiased estimators
are established. A geometric characterization of a distribution that optimizes the performance of these estimators is given. An algorithm based on convex hulls that finds the optimal  distribution truncated to its first \(m\) terms in time linear in \(m\) is described.  The algorithm exploits a connection with a  recent randomized dimension reduction method and is illustrated via a numerical example.
     }
\commentt{
}{\maketitle}

\commentt{}
{Subject classifications:  Simulation: Efficiency; Analysis of algorithms: Computational complexity;
Area of Review: Simulation}
\section{Introduction}
Monte Carlo methods are used in a variety of  domains  such as financial engineering, queuing  networks,  and machine learning. In general, however, Monte Carlo methods are computationally costly. Variance reduction techniques such as importance sampling, control variate methods, stratified sampling and splitting techniques can significantly  improve the efficiency of Monte Carlo methods (e.g.~\cite{glasserman2004Monte,asmussenGlynn2007,rubinstein2016simulation}). The multilevel Monte Carlo 
 method (MLMC), introduced by~\citet{Giles2008},   dramatically reduces the computational cost of estimating an expected value arising
from a stochastic differential equation. \citet{mcleish2011}, \citet{glynn2014exact} and    \citet{GlynnRhee2015unbiased}    provide related randomized multilevel Monte Carlo methods (RMLMC) that produce  unbiased estimators" for  equilibrium expectations of functionals on Markov chains, and for expectations of functionals  arising in    stochastic differential equations. Sufficient conditions guaranteeing the  square-integrability of these estimators are given in \cite{mcleish2011,GlynnRhee2015unbiased}. \citet{jacob2015nonnegative} study the existence of unbiased nonnegative estimators.       RMLMC  and related methods  have been used in a variety of  contexts such as  the unbiased estimation of a function of the mean of a random variable~\cite{blanchet2015unbiased,kroese2019unbiased}, the design of Markov chain Monte Carlo methods~\cite{bardenet2017,agapiou2018unbiased,middleton2018unbiased}, unbiased inference for hidden Markov models~\cite{vihola2018unbiasedInference}, pricing of Asian options under general models~\cite{kahale2018Asian}, and stochastic optimization~\cite{blanchet2019unbiased}.  \citet{Vihola2018} describes  stratified  RMLMC methods that, under certain conditions, are shown to be asymptotically as
efficient as MLMC. The randomized dimension reduction method,   recently  introduced in \cite{kahNIPS16,kahaRandomizedDimensionReduction19}, is another technique that can provably achieve substantial  variance reduction in high-dimensional settings,   such as the estimation of the expectation of a functional of a time-varying Markov chain at a long horizon. 

 Using the terminology of \citet{GlynnRhee2015unbiased}, the ``coupled sum'' and ``independent sum'' unbiased estimators take as parameter the distribution of an integral random variable. 
In   \cite[Section 3]{GlynnRhee2015unbiased}, an algorithm that finds in \(O(m^3)\) time  an \(m\)-truncated distribution that optimizes the efficiency of these estimators is given. 
On the other hand, the asymptotic efficiency of the randomized dimension reduction method is maximized  in \cite{kahaRandomizedDimensionReduction19} via a new geometric algorithm that   solves an \(m\)-dimensional optimization problem in \(O(m)\) time.    \citet{kahaRandomizedDimensionReduction19} points out that the same geometric algorithm  solves the optimization
problem  in \cite[Section 3]{GlynnRhee2015unbiased}      in \(O(m)\) time.

The output of unbiased estimators can be analysed using  well-known tools.
For instance, drawing independent copies of an unbiased estimator allows the construction of normal confidence intervals \cite[Section III.1]{asmussenGlynn2007}.  In addition, this construction is easy to parallelize. Also, \citet{glynn1992asymptotic} have established a  central limit theorem on the average of independent copies of an unbiased estimator under a  computational budget constraint.   On the other hand, the output of biased estimators can be difficult to analyse, even for estimators which are asymptotically consistent (e.g. \cite[Section 4.5.1]{glasserman2004Monte}). This highlights the importance of the RMLMC techniques. Motivated
by these considerations and by the wide range of applications of the RMLMC methods, this note studies the coupled sum  and independent sum estimators  in a general framework and makes three main contributions:
\begin{enumerate}
\item 
It gives a necessary and sufficient condition for the square-integrability of  the coupled sum  (resp. independent sum) estimator.
When this condition is met, it is shown that the corresponding estimator is unbiased and an expression for the second moment is derived.  
An example showing that the new conditions are strictly weaker than the sufficient condition  in    \cite{GlynnRhee2015unbiased} is given.

  \item 
Under general conditions,  it gives a simple geometric characterization, based on convex hulls,  of  distributions with infinite support that optimize the performance of the coupled sum  (resp. independent sum)  estimator.     \citet[Theorem 3]{GlynnRhee2015unbiased} show that such a distribution can be found by solving a certain combinatorial problem,
but do not provide the solution to this problem in the infinite support case. \item
 Building on techniques developed in \cite{kahaRandomizedDimensionReduction19}, it describes an algorithm that   finds an optimal \(m\)-truncated  distribution for each of  these estimators in \(O(m)\)  time. The algorithm, based on convex hulls, is simple to implement.   \citet{GlynnRhee2015unbiased} give an alternative algorithm       based on dynamic programming that runs in \(O(m^3)\) time. More recently, \citet{cui2019optimal} give yet another algorithm that solves this problem in \(O(m)\)  time by using a dual formulation of the optimization problem.  
\end{enumerate}    
The rest of the paper is organised as follows. \S\ref{se:integrabilityConditions} gives necessary and sufficient conditions for the square-integrability of the coupled sum  and independent sum estimators, and presents expressions for the second moment of these estimators when these conditions are met.
\S\ref{se:optimalInfSupport} describes a geometric characterization of distributions with infinite support that optimize the efficiency of these estimators.
\S\ref{se:optimalFiniteSupport} shows how to calculate in \(O(m)\) time an optimal \(m\)-truncated distribution. \S\ref{se:numer} describes a numerical example. \S\ref{se:conc} contains concluding remarks.
Omitted proofs are in the appendix.

\section{The coupled and independent sum estimators}\label{se:integrabilityConditions}   Let \(\mathbb{R}_{+}\)  denote the set of nonnegative real numbers.  For  a square-integrable random variable  \(X\), let \(||X||=\sqrt{E(X^{2})}\). The  coupled and independent sum estimators efficiently estimate  the expectation of a random variable \(Y\)   that is approximated by random variables \(Y_{n}\),  \(n\geq0\). By convention,    \(Y_{-1}=0\).  Let \((\tilde\Delta_{n}:n\ge0)\) be a sequence of independent random variables such that \(\tilde\Delta_{n}\) has the same distribution as \(Y_{n}-Y_{n-1}\)  for \(n\ge0\).
It is assumed throughout the paper that  \(Y\) and  \(Y_{n}\) are   square-integrable, that \(||Y_{n}-Y||\)   goes to \(0\) as \(n\) goes to infinity, and that
the expected time to generate \((Y_{0},\dots,Y_{n})\) (resp.  \((\tilde\Delta_{0},\dots,\tilde\Delta_{n})\)) is finite, for \(n\geq0\). Let  \begin{displaymath}
A=\{(q_{i}: i\ge0)\in\mathbb{R}^{\mathbb{N}}:q_{0}=1,  q_{i}\ge  q_{i+1}>0 \text{ for }i\ge0,\lim_{i\rightarrow\infty} q_{i}=0\}.
\end{displaymath}
Example~\ref{ex:option} is a standard application of MLMC and RMLMC methods.
\begin{example}\label{ex:option}
Let \(T\) be a fixed maturity and let  \((X(t):0\leq t\leq T)\) be a stochastic process that solves the stochastic differential equation
\begin{equation*}
dX(t)= a(X(t),t)\,dt+b(X(t),t)\,dW,
\end{equation*} where \(a\) and \(b\) are real-valued functions on \(\mathbb{R}\times[0,T]\)  and \(W\) is a one-dimensional Brownian motion. Suppose that, for \((x,t)\in \mathbb{R}\times[0,T]\), \(a(x,t)\) and 
 \(b(x,t)\) can be calculated in constant time.   Option pricing applications often need to estimate \(E(Y)\), where  \(Y=f(X(T))\) and \(f\) is a payoff function. The process \(X\) can be approximately simulated via the Milstein   discretisation scheme as follows. For  \(i\ge0\), define  recursively the sequence \((X^{(i)}_{k}:0\leq k\leq 2^{i})\)   by setting \(X^{(i)}_{0}=X(0)\) and
\begin{equation*}
 X^{(i)}_{k+1}=X^{(i)}_{k}+a(X^{(i)}_{k},k\Delta t)\,\Delta t+b(X^{(i)}_{k},k\Delta t)\,\Delta W+\frac{1}{2}b(X^{(i)}_{k},k\Delta t)\frac{\partial b}{\partial x}(X^{(i)}_{k},k\Delta t)((\Delta W)^{2}-\Delta t),
\end{equation*} 
\(0\leq k\leq 2^{i}-1\), where  \(\Delta t=2^{-i}T\) and \(\Delta W=W((k+1)\Delta t)-W(k\Delta t)\). Let \(Y_{i}=f(X^{(i)}_{2^{i}})\). Under certain conditions on \(a\), \(b\) and \(f\), \begin{equation}\label{eq:exOptionYiY}
||Y_{i}-Y||^{2}\leq c2^{-2i},
\end{equation}where \(c\) is a constant (e.g. \citep{kloedenPlaten1992}). Note that \(Y_{i}\) can be calculated in \(O(2^{i})\) time.
\end{example}   
\subsection{The coupled sum estimator}\label{sub:coupledSum}   Let  \(q\in A\)  and  let \(N\)   be an integral random variable     independent of  \((Y_{i}: i\geq0)\)    such that \(\Pr(N\ge i)= q_{i}\)   for \(i\ge0\). Following \citet{mcleish2011} and~\citet{GlynnRhee2015unbiased}, define the coupled sum estimator as
\begin{equation*}
 \bar Z = \sum^{N}_{i=0}\frac{Y_{i}-Y_{i-1}}{ q_{i}}.
\end{equation*}   
Theorem~\ref{th:GenCase}  gives a necessary and sufficient condition for   \(\bar Z\) to be square-integrable.  Moreover, under this condition, it shows that   \(\bar Z\) is an unbiased estimator for \(E(Y)\), and gives an expression for the second moment of   \(\bar Z\). Note that each term in the LHS of   \eqref{eq:NewCondZIntegral} is nonnegative because   \(q\in A\).

  \begin{theorem}\label{th:GenCase}
 The coupled sum estimator  \(\bar Z\) is square-integrable if and only if
\begin{equation}\label{eq:NewCondZIntegral}
\sum^{\infty}_{i=0}(\frac{1}{q_{i+1}}-\frac{1}{q_{i}})||Y_{i}-Y||^{2}<\infty.
\end{equation}
Furthermore, if   \eqref{eq:NewCondZIntegral} holds then \(E(\bar Z)=E(Y)\) and
  \begin{equation}\label{eq:normbarZ}
||\bar Z||^{2}=||Y||^{2}+\sum^{\infty}_{i=0}(\frac{1}{q_{i+1}}-\frac{1}{q_{i}})||Y_{i}-Y||^{2}.
\end{equation}
\end{theorem}  
Theorem 1 in~\cite{GlynnRhee2015unbiased} gives a sufficient condition  for   \(\bar Z\) to be square-integrable. More precisely, it shows that if \begin{equation}
\label{eq:condThGlynn} \sum^{\infty}_{i=1}\frac{||Y_{i-1}-Y||^{2}}{ q_i}<\infty,
\end{equation}then     \(\bar Z\) is square-integrable,  and   \(\bar Z\) is an unbiased estimator for \(E(Y)\), and gives an expression for  \(||\bar Z||^{2}\).    
However, the following toy example shows that    \eqref{eq:condThGlynn} is not a necessary condition for \(\bar Z\) to be square-integrable.
\begin{example}
\label{ex:toy}Assume that \(Y=0\) and \(Y_{i}=(i+1)^{-3/2}\), and
 \(q_{i}=(i+1)^{-2}\)
for \(i\ge0\). Then  \eqref{eq:NewCondZIntegral} holds and so   \(\bar Z\) is square-integrable, but      \eqref{eq:condThGlynn} does not hold.
\end{example}
\citet{mcleish2011} gives alternative conditions that guarantee the unbiasedness  and  square-integrability of    \(\bar Z\), and provides an alternative expression for   \(||\bar Z||^{2}\).   
\subsection{The independent sum estimator}
Let  \(q\in A\)  and  let \(N\)   be an integral random variable     independent  of \((\tilde\Delta_{i}: i\geq0)\)  such that \(\Pr(N\ge i)= q_{i}\)   for \(i\ge0\). Following~\citet{GlynnRhee2015unbiased}, define the independent sum estimator as
\begin{equation*}
 \tilde Z = \sum^{N}_{i=0}\frac{\tilde\Delta_{i}}{ q_i}.
\end{equation*}
 Theorem~\ref{th:GenCaseTildeZ} below gives a necessary and sufficient condition for   \(\tilde Z\) to be square-integrable.  Furthermore, if this condition is met, it shows that   \(\tilde Z\) is an unbiased estimator for \(E(Y)\), and gives an expression for the second moment of   \(\tilde Z\).
\begin{theorem}\label{th:GenCaseTildeZ}
The independent sum estimator  \(\tilde Z\) is square-integrable if and only if
\begin{equation}\label{eq:NewCondZIntegralTilde}
\sum^{\infty}_{i=0}\left(\frac{\var(Y_{i}-Y_{i-1})}{q_{i}}+(\frac{1}{q_{i+1}}-\frac{1}{q_{i}})(E(Y_{i}-Y))^{2}\right)<\infty.
\end{equation}
Furthermore, if \eqref{eq:NewCondZIntegralTilde} holds then \(E(\tilde Z)=E(Y)\) and
  \begin{equation}\label{eq:normTildeZ}
||\tilde Z||^{2}=(E(Y))^{2}
+\sum^{ \infty}_{i=0}\left(\frac{\var(Y_{i}-Y_{i-1})}{q_{i}}
+(\frac{1}{q_{i+1}}-\frac{1}{q_{i}})(E(Y_{i}-Y))^{2}\right).
\end{equation}
\end{theorem}
Theorem~2 in~\cite{GlynnRhee2015unbiased} shows that    \eqref{eq:condThGlynn}  is a sufficient condition  for the square-integrability of   \(\tilde Z\). Assuming that    \eqref{eq:condThGlynn} holds, Theorem~2 in~\cite{GlynnRhee2015unbiased}    shows that     \(\tilde Z\) is square-integrable,  and   \(\tilde Z\) is an unbiased estimator for \(E(Y)\), and gives an expression for \(||\tilde Z||^{2}\). 
 In Example~\ref{ex:toy}, however,   \eqref{eq:NewCondZIntegralTilde} holds and so  \(\tilde Z\) is square-integrable, but     \eqref{eq:condThGlynn}    does not hold. Thus    \eqref{eq:condThGlynn} is not a necessary condition for \(\tilde Z\) to be square-integrable.\section{Optimal distribution: the infinite support case}\label{se:optimalInfSupport}
\subsection{The coupled sum estimator}\label{sub:optCoupledSum}
For \(i\geq 0\), let \(\bar t_{i}\) be the expected cost required to simulate the sequence \((Y_{0},\dots,Y_{i-1})\). By convention, \(\bar t_{0}=0\). 
    It is assumed that the sequence \(\bar t=(\bar t_{i}:i\ge0)\) is strictly increasing and that \(\bar t_{i}\)
  goes to infinity as \(i\) goes to infinity. Let   \( \bar\tau\)  be the time required to generate     \(\bar Z\).  As observed by~\citet{GlynnRhee2015unbiased},
\begin{equation*}
E(\bar\tau)=E(\bar t_{N+1})=\sum^{\infty}_{i=0} q_{i} (\bar t_{i+1}-\bar t_{i}).
\end{equation*}Note that \(E(\bar\tau)\)  is either   infinite  or is finite and positive.
In Example~\ref{ex:option}, assuming that \eqref{eq:exOptionYiY} holds,
  \(q\) can be chosen so that  \(\bar Z\) is square-integrable  and  \(E(\bar\tau)\) is  finite by setting \(q_{i}=2^{-3i/2}\). Indeed,  for this choice of \(q\),   \eqref{eq:NewCondZIntegral} holds and   \(E(\bar\tau)\) is  finite since \(\bar t_{i}=O(2^{i})\).   
   
 Glynn and Whitt (1992) show that the efficiency of an unbiased estimator is inversely
proportional to the product of the variance and expected running
time. Thus, maximizing the efficiency of  \(\bar Z\) amounts to finding a sequence \(q\in A\) that satisfies \eqref{eq:NewCondZIntegral} and  minimizes \(E(\bar\tau)\var(\bar Z)\). 
For   \(q\in A\), and any strictly increasing sequence \(\vartheta=(\vartheta_{0},\vartheta_{1},\dots)\in\{0\}\times\mathbb{R}_{+}^{\mathbb{N}}\), and \(\gamma=(\gamma_{0},\gamma_{1},\dots)\in\mathbb{R}_{+}^{\mathbb{N}}\),   set
\begin{equation}\label{eq:RDefGen}
R(q;\vartheta,\gamma)=\lim_{n\rightarrow\infty}(\sum^{n}_{i=0}q_{i}(\vartheta_{i+1}-\vartheta_{i}))(\gamma_{0}+\sum^{n}_{i=1}(\frac{1}{q_{i}}-\frac{1}{q_{i-1}})\gamma_{i}).
\end{equation}
Thus, \(R(q;\vartheta,\gamma)\) is the limit of a nonnegative increasing sequence, and so  \(R(q;\vartheta,\gamma)\) is   either infinite or is nonnegative and finite.  

Define the sequence \(\bar\mu=(\bar\mu_{i}: i\geq0)\), where  \(\bar\mu_0=\var(Y)\) and \(\bar\mu_{i}=||Y_{i-1}-Y||^{2}\) for \(i\ge1\).  For simplicity, this subsection assumes  that the sequence \(\bar\mu\) is positive. The case where \(Y=Y_{m}\) for some integer \(m\) is studied in \S\ref{sub:trcoupledsum}.  By Theorem~\ref{th:GenCase}, if   \(q\in A\) is such that \(\bar Z\) is square-integrable  and  \(E(\bar\tau)\) is finite, then \begin{equation}\label{eq:varBarZ}
\var(\bar Z)=\bar\mu_0+\sum^{\infty}_{i=1}(\frac{1}{q_{i}}-\frac{1}{q_{i-1}})\bar\mu_{i},
\end{equation} 
 and so, by  \eqref{eq:RDefGen},    \(R(q;\bar t,\bar\mu)=E(\bar\tau)\var(\bar Z)\) is finite. Conversely, if   \(q\in A\) and \(R(q;\bar t,\bar\mu)\) is finite, then by  \eqref{eq:RDefGen}, \(E(\bar\tau)\) is finite and   \eqref{eq:NewCondZIntegral} holds. Hence  \(\bar Z\) is square-integrable  and   \(R(q;\bar t,\bar\mu)=E(\bar\tau)\var(\bar Z)\). Thus, maximizing the efficiency of \(\bar Z\) amounts to finding a sequence \(q\in A\) that minimizes \(R(q;\bar t,\bar\mu)\). Under general conditions, given a strictly increasing sequence \(\vartheta=(\vartheta_{0},\vartheta_{1},\dots)\in\{0\}\times\mathbb{R}_{+}^{\mathbb{N}}\) and  a sequence    \(\gamma\in\mathbb{R}_{+}^{\mathbb{N}}\),  Theorem~\ref{th:OptTimeVarianceGen} below gives a characterization of a sequence \(q^{*}\) that minimizes \(R(q;\vartheta,\gamma)\) under the constraint that \(q\in A\). The sequence \(q^{*}\) clearly depends on \(\vartheta\) and \(\gamma\).
Theorem~\ref{th:OptTimeVarianceGen}  generalizes Theorem~3 in~\cite{kahaRandomizedDimensionReduction19} to infinite sequences.  

Let \(\gamma'\in\mathbb{R}_{+}^{\mathbb{N}}\) be such that the set  \(\{(\vartheta_{i},\gamma'_{i}):i\in
\mathbb{N}\}\)
forms the lower hull of the set  \(\{(\vartheta_{i},\gamma_{i}):i\in
\mathbb{N}\}\). Thus  \(\gamma'\) is the supremum  of all real sequences  such that \(\gamma'\le\gamma\) and  the sequence   \((\theta_{i})\) is  increasing, where  
\begin{equation}\label{eq:defthetaGen}
\theta_{i}=\frac{\gamma'_{i+1}-\gamma'_{i}}{\vartheta_{i+1}-\vartheta_{i}},
\end{equation}\( i\geq0\).   
 
\begin{theorem}\label{th:OptTimeVarianceGen}
Let  \(\vartheta=(\vartheta_{0},\vartheta_{1},\dots)\in\{0\}\times\mathbb{R}_{+}^{\mathbb{N}}\) be a  strictly increasing sequence, and let  \(\gamma=(\gamma_{0},\gamma_{1},\dots)\in\mathbb{R}_{+}^{\mathbb{N}}\)     be a positive sequence that
goes to \(0\) at infinity.  Assume there is \(q\in A\) such that \(R(q;\vartheta,\gamma)\) is finite. For \(i\geq 0\),  set  \(q^{*}_{i}=\sqrt{\theta_{i}/\theta_{0}}\), where \(\theta_{i}\) is given by~\eqref{eq:defthetaGen}. Then \(q^{*}=(q^{*}_{i}: i\ge0)\in A\) and    \(R(q^{*};\vartheta,\gamma)\) is finite. Moreover  \(q^{*}=\arg\min_{q\in A}R(q;\vartheta,\gamma)\)   and 
\begin{equation}\label{eq:optTimeVariance}
 R(q^{*};\vartheta,\gamma)=\bigg(\sum_{i=0}^{\infty}\sqrt{(\gamma'_{i}-\gamma'_{i+1})(\vartheta_{i+1}-\vartheta_{i})}\bigg
)^2,\end{equation} where the series in the RHS of \eqref{eq:optTimeVariance} is convergent.
\end{theorem}
Thus, if there is \(q\in A\) such that \(\bar Z\) is square-integrable  and   \(E(\bar\tau)\) is finite,   then the efficiency of \(\bar Z\) is maximized when \(\Pr(N\ge i)= q^{*}_{i}\)   for \(i\ge0\), where \(q^{*}\) is the sequence described in  Theorem~\ref{th:OptTimeVarianceGen} 
 with \(\vartheta=\bar t\) and \(\gamma=\bar\mu\). 
\subsection{The independent sum estimator}For \(i\geq 0\), let \(\tilde t_{i}\) be the expected cost required to simulate the sequence \((\tilde\Delta_{0},\dots,\tilde\Delta_{i-1})\). By convention, \(\tilde t_{0}=0\). 
    Assume that the sequence \(\tilde t=(\tilde t_{i}:i\ge0)\) is strictly increasing and that \(\tilde t_{i}\)
  goes to infinity as \(i\) goes to infinity. Let   \( \tilde\tau\)  be the time required to generate     \(\tilde Z\).  Then\begin{equation*}
E(\tilde\tau)=E(\tilde t_{N+1})=\sum^{\infty}_{i=0} q_{i} (\tilde t_{i+1}-\tilde t_{i}).
\end{equation*}
\comment{
\begin{example}
Assume that there are constants \(c\) and \(c'\) and \(\delta<-2\) such that \(\tilde t_{i}\leq ci\) and \(||Y_{i}-Y_{i-1}||^{2}\le c'(i+1)^{\delta}\) for \(i\ge0\). Let \(q_{i}=(i+1)^{\delta/2}\) for \(i\geq0\). By \eqref{eq:tildeTau}, \begin{equation}
E(\tilde\tau)\le cE(N+1)=c\sum^{\infty}_{i=0} q_{i}
\end{equation}is finite.
By the Cauchy-Schwartz inequality,
for \(0\le i\leq j\), \begin{displaymath}
(E(Y_{j}-Y_{i}))^{2}\le\left(\sum^{j}_{k=i+1}k^{\delta/2}\right)\left(\sum^{j}_{k=i+1}\frac{(E(Y_{k}-Y_{k-1}))^{2}}{k^{\delta/2}}\right).
\end{displaymath}
Thus, \eqref{eq:NewCondZIntegral} holds and, by Theorem~\ref{th:GenCase}, \(\bar Z\) is square-integrable. \end{example}
}
In Example~\ref{ex:option}, assuming that \eqref{eq:exOptionYiY} holds,
  \(q\) can be chosen so that  \(\tilde Z\) is square-integrable  and  \(E(\tilde\tau)\) is  finite by setting \(q_{i}=2^{-3i/2}\). Indeed, for \(i\ge0\),\begin{displaymath}
|E(Y_{i}-Y)|\leq||Y_{i}-Y||\leq\sqrt{c}2^{-i}
\end{displaymath}
and, for \(i\ge1\),\begin{displaymath}
\std(Y_{i}-Y_{i-1})\leq||Y_{i}-Y||+||Y_{i-1}-Y||\leq 3\sqrt{c}2^{-i}.
\end{displaymath} Thus, for this choice of \(q\),   \eqref{eq:NewCondZIntegralTilde} holds and   \(E(\tilde\tau)\) is  finite since \(\tilde t_{i}=O(2^{i})\).
 
When \(\tilde Z\) is square-integrable, Proposition~\ref{pr:varZtilde}  gives an expression for  \(\var(\tilde Z)\) similar to that of \(\var(\bar Z)\) in~\eqref{eq:varBarZ}.

\begin{proposition}\label{pr:varZtilde}If \(q\in A\) is such that \(\tilde Z\) is square-integrable, then   \(\tilde\mu_0=\sum^{\infty}_{j=0}\var(Y_{j}-Y_{j-1})\) is finite, and  
\begin{equation*}
\var(\tilde Z)=\tilde\mu_0+\sum^{\infty}_{i=1}(\frac{1}{q_{i}}-\frac{1}{q_{i-1}})\tilde\mu_{i},
\end{equation*} 
where, 
for \(i\ge1\),
\begin{displaymath}
\tilde\mu_{i}=(E(Y_{i-1}-Y))^{2}+\sum^{\infty}_{j=i}\var(Y_{j}-Y_{j-1}).
\end{displaymath}\end{proposition} 
Maximizing the efficiency of  \(\tilde Z\) amounts to finding a sequence \(q\in A\) that satisfies   \eqref{eq:NewCondZIntegralTilde} and  minimizes \(E(\tilde\tau)\var(\tilde Z)\). Suppose  that  \(q\in A\)  is such that \(\tilde Z\) is square-integrable  and  \(E(\tilde\tau)\) is finite. By Proposition~\ref{pr:varZtilde}, and since \(|E(Y_{i-1}-Y)|\leq ||Y_{i-1}-Y||\), the sequence  \(\tilde\mu=(\tilde\mu_{i}: i\geq0) \)  goes to \(0\) at infinity. For simplicity, assume that   \(\tilde\mu\)   is positive.   Using arguments similar to those used in \S\ref{sub:optCoupledSum}, it can be shown that the efficiency of \(\tilde Z\) is maximized when \(\Pr(N\ge i)= q^{*}_{i}\)   for \(i\ge0\), where \(q^{*}\) is the sequence described in  Theorem~\ref{th:OptTimeVarianceGen} 
 with \(\vartheta=\tilde t\) and \(\gamma=\tilde\mu\). 

  \section{Optimal distribution: the truncated case}\label{se:optimalFiniteSupport}
    \subsection{A truncated coupled sum estimator} \label{sub:trcoupledsum}
Fix \(m\geq1\) and let  \begin{displaymath}
A^{(m)}=\{(q_{0},\dots,q_{m})\in\mathbb{R}^{m+1}:1=q_{0}\ge q_{1}\ge\cdots\ge q_{m}>0\}.
\end{displaymath}Let  \(q\in A^{(m)}\)  and  let \(S\)   be an integral random variable  in \(\{0,\dots,m\}\)   independent of  \((Y_{i}: 0\leq i\leq m)\)  such that \(\Pr(S\ge i)= q_{i}\)   for \(0\le i\le m\).
Set
\begin{equation*}
 \bar Z^{(m)} = \sum^{S}_{i=0}\frac{Y_{i}-Y_{i-1}}{q_{i}}.
\end{equation*}Let  \(N\) be an integral random variable with infinite support independent of  \((Y_{i}: 0\leq i\leq m)\)  such that   \(\Pr(N\ge i)=q_{i}\)   for \(0\le i\le m\). Thus, \(S\) has the same distribution as \(N\wedge m\). By applying Theorem~\ref{th:GenCase} to \(N\) and to the sequence \((Y_{n\wedge m}: n\geq0)\), with \(Y=Y_{m}\),  it follows that  \( \bar Z^{(m)}\) is square-integrable, with \(E( \bar Z^{(m)})=E(Y_{m})\) and
  \begin{eqnarray}\label{eq:truncated}\nonumber  
||\bar Z^{(m)}||^{2}&=&  
||Y_{m}||^{2}+\sum^{m-1}_{i=0}(\frac{1}{q_{i+1}}-\frac{1}{q_{i}})||Y_{i}-Y_{m}||^{2}
\\&=&\sum^{m}_{i=0}\frac{||Y_{i-1}-Y_{m}||^{2}-||Y_{i}-Y_{m}||^{2}}{ q_{i}}.
\end{eqnarray} 
Define the sequence  \(\bar\eta=(\bar\eta_i:  0\le i\le m+1)\), where \(\bar\eta_0=\var(Y_{m})\) and \(\bar\eta_{i}=||Y_{i-1}-Y_{m}||^{2}\) for \(1\le i\le m+1\). Thus, \(\bar\eta\) can be considered as the truncated counterpart of the sequence 
\(\bar \mu\) defined in \S\ref{sub:optCoupledSum}.   Note that \(\bar\eta_{m+1}=0\). Assume for simplicity that \(\bar\eta_{i}>0\) for  \(0\leq i\leq m\). 
By \eqref{eq:truncated}, 
\begin{equation*}
\var( \bar Z^{(m)})=\sum^{m}_{i=0}\frac{\bar\eta_{i}-\bar\eta_{i+1}}{ q_{i}}.
\end{equation*}
Let \(\bar\tau^{(m)}\) be the expected time to simulate \(\bar Z^{(m)}\).  As \(\Pr(S=i)=q_{i}-q_{i+1}\) for \(0\leq i\leq m\), where \(q_{m+1}=0\) by convention, \begin{equation*}
E(\bar\tau^{(m)})=\sum^{m}_{i=0}(q_{i}-q_{i+1})\bar t_{i+1}=\sum^{m}_{i=0}q_{i}(\bar t_{i+1}-\bar t_{i}).
\end{equation*}\subsection{The optimal distribution}

Maximizing the efficiency of  \(\bar Z^{(m)}\) amounts to finding a sequence \(q\in A^{(m)}\) that   minimizes \(E(\bar\tau^{(m)})\var(\bar Z^{(m)})\).
For \(q\in A^{(m)}\), and any strictly increasing sequence  \(\vartheta=(\vartheta_{0},\dots,\vartheta_{m+1})\in\{0\}\times \mathbb{R}_{+}^{m+1}\),  and \(\gamma=(\gamma_{0},\dots,\gamma_{m+1})\in\mathbb{R}_{+}^{m+1}\times\{0\}\),     set  
   \begin{equation*}
R^{(m)}(q;\vartheta,\gamma)=(\sum^{m}_{i=0}q_{i}(\vartheta_{i+1}-\vartheta_{i}))(\sum^{m}_{i=0}\frac{\gamma _{i}-\gamma _{i+1}}{q_{i}}).
\end{equation*} 
Hence, \(E(\bar\tau^{(m)})\var(\bar Z^{(m)})=R^{(m)}(q;\bar t^{(m)},\bar\eta)\), where \(\bar t^{(m)}=(\bar t_{0},\dots,\bar t_{m+1})\). Thus, we need to calculate a sequence \(q\in A^{(m)}\) that minimizes \(R^{(m)}(q; \bar t^{(m)},\bar\eta)\). Given a sequence    \(\gamma\in\mathbb{R}_{+}^{m+1}\times\{0\}\) whose first \(m+1\) components are positive,  Theorem~\ref{th:OptTimeVariance}, a direct consequence of Theorem~3 in \cite{kahaRandomizedDimensionReduction19}, shows how to calculate in \(O(m)\) time a vector \(q^{*}\) that minimizes \(R^{(m)}(q;\bar t^{(m)},\gamma)\) under the constraint that \(q\in A^{(m)}\). The vector \(q^{*}\) clearly depends on \(\bar t^{(m)}\) and on \(\gamma\).

Let \(\gamma'=(\gamma'_{0},\dots,\gamma'_{m+1})\in\mathbb{R}^{m+2}\) be such that the set  \(\{(\bar t_{i},\gamma'_{i}):0\leq i\leq m+1\}\)
forms the lower hull of the set  \(\{(\bar t_{i},\gamma_{i}):0\leq i\leq m+1\}\). In other words,  \(\gamma'\) is the supremum  of all sequences in \(\mathbb{R}^{m+2}\) such that \(\gamma'\le\gamma\) and  the sequence   \((\theta_{i})\) is  increasing, where  
\begin{equation}\label{eq:deftheta}
\theta_{i}=\frac{\gamma'_{i+1}-\gamma'_{i}}{\bar t_{i+1}-\bar t_{i}},
\end{equation}  \(0\leq i\leq m\).  For instance, if \(m=5\), with \(\bar t_{i}=i\) and \(\gamma=(20,22,14,5,4,1,0)\), then  \(\gamma'=(20,15,10,5,3,1,0)\), as illustrated in Fig.~\ref{fig:hull}.
 The vector  \(\gamma'\) can be calculated in \(O(m)\) time   using the convex hull algorithm of \citet{Andrew79}. 

\begin{figure}
\caption{Lower hull}
\begin{center}
\begin{tikzpicture}
\begin{axis}[
    title={},
    xlabel={$i$},
    xmin=0, xmax=6,
    ymin=0, ymax=30,
    xtick={0,1,2,3,4,5,6},
    ytick={0,5,10,15,20,25},
    legend pos=north east,
    ymajorgrids=true,
    grid style=dashed,
]
 
\addplot[
    color=blue,
    mark=square,
    ]
    coordinates {
    (0,20)(1,22)(2,14)(3,5)(4,4)(5,1)(6,0)
    };
 \addplot[
    color=red,
    mark=star,
    ]
    coordinates {
    (0,20)(1,15)(2,10)(3,5)(4,3)(5,1)(6,0)
    };
    \legend{$\gamma$,$\gamma'$}

\end{axis}
\end{tikzpicture}
\end{center}
\label{fig:hull}
\end{figure}
\begin{theorem}[\cite{kahaRandomizedDimensionReduction19}]\label{th:OptTimeVariance}
Let   \(\gamma\)    be a vector in \(\mathbb{R}^{m+1}\times\{0\}\) whose first \(m+1\) components are positive. For \(0\leq i\leq m\),  set  \(q^{*}_{i}=\sqrt{\theta_{i}/\theta_{0}}\), where \(\theta_{i}\) is given by~\eqref{eq:deftheta}, and let \(q^{*}=(q^{*}_{0},\dots,q^{*}_{m})\). Then \(q^{*}=\arg\min_{q\in A^{(m)}}R^{(m)}(q;\bar t^{(m)},\gamma)\).
\end{theorem}
The proof of  \cite[Theorem 3]{kahaRandomizedDimensionReduction19} shows that the conclusions of Theorem~\ref{th:OptTimeVariance} are valid for any strictly increasing sequence \(\bar t^{(m)}\in\{0\}\times\mathbb{R}_{+}^{m+1}\). 
\subsection{The algorithm description}\label{sub:hull}
   Combining the previously discussed elements yields  an algorithm that takes as input  the vectors  \(\bar t^{(m)}\) and  \(\bar\eta\), and outputs \(q^{*}=\arg\min_{q\in A_{m}}R^{(m)}(q;\bar t^{(m)},\bar\eta)\) in \(O(m)\) time. The first two steps of the algorithm are adapted from \citet{Andrew79}.  The  algorithm  first generates by backward induction a subset \(B(j)\) of  \(\{j,\dots,m+1\}\), \(0\leq j\leq m\), so that the lower hull  \(\{(\bar t_{i},\bar\eta'_{i}):0\leq i\leq m+1\}\)
of the set  \(\{(\bar t_{i},\bar\eta_{i}):0\leq i\leq m+1\}\) is  obtained by piece-wise interpolation of \(\bar\eta\) on  \(B(0)\). More precisely,   \(B(0)\) contains   \(\{0,m+1\}\), and if  \(i'\) and \(i''\) are two consecutive  elements of \(B(0)\) with \(i'\leq i\le i''\), then  \((\bar t_{i},\bar\eta'_{i})\) lies on the segment \([(\bar t_{i'},\bar\eta_{i'}),(\bar t_{i''},\bar\eta_{i''})]\). The third step calculates \(q^{*}\) via \(B(0)\). \begin{enumerate}
\item 
Set \(B(m)=\{m,m+1\}\).
\item For \(j=m-1\) down to \(0\), denote by \(i_{0}<\cdots<i_{l}\)  the elements of \(B(j+1)\).
Determine the smallest  element \(k\) of the set   \(\{0,\dots,l-1\}\)   such that \((\bar t_{i_k},\bar\eta_{i_{k}})\) lies below the segment
\([(\bar t_{j},\bar\eta_{j}),(\bar t_{i_{k+1}},\bar\eta_{i_{k+1}})]\), if such \(k\) exists, otherwise let \(k=l\).  Set  \(B(j)=\{j,i_{k},\dots,i_{l}\}\). 
\item
For \(i=0\) to \(m\), let \(i'\) and \(i''\) be two consecutive  elements of \(B(0)\) with \(i'\leq i< i''\). Set  
\begin{equation*}
\theta_{i}=\frac{\gamma_{i''}-\gamma_{i'}}{\bar t_{i''}-\bar t_{i'}}
\end{equation*} and   \(q^{*}_{i}=\sqrt{\theta_{i}/\theta_{0}}\). 
\end{enumerate}
As pointed out in \cite{GlynnRhee2015unbiased,Vihola2018}, the vector   \(\bar\eta\) is not known exactly, in general, but can be estimated by Monte Carlo simulation.
\subsection{A truncated independent sum estimator}\label{sub:trIndSum}
Let  \(q\in A^{(m)}\)  and  let \(S\)   be an integral random variable  in \([0,m]\)   independent of  \((\tilde\Delta_{i}: 0\leq i\leq m)\) such that \(\Pr(S\ge i)= q_{i}\)   for \(0\le i\le m\). Set
\begin{equation*}
 \tilde Z^{(m)} = \sum^{S}_{i=0}\frac{\tilde \Delta_{i}}{q_{i}}.
\end{equation*}By Theorem~\ref{th:GenCaseTildeZ}, Proposition~\ref{pr:varZtilde}, and arguments  similar to those used in \S\ref{sub:trcoupledsum},   it follows that  \( \tilde Z^{(m)}\) is square-integrable, with \(E( \tilde Z^{(m)})=E(Y_{m})\), and
\begin{equation*}
\var( \tilde Z^{(m)})=\sum^{m}_{i=0}\frac{\tilde\eta_{i}-\tilde\eta_{i+1}}{ q_{i}},
\end{equation*}
where \(\tilde\eta_0=\sum^{m}_{j=0}\var(Y_{j}-Y_{j-1})\) and, 
for \(1\le i\le m+1\),
\begin{displaymath}
\tilde\eta_{i}=(E(Y_{i-1}-Y_{m}))^{2}+\sum^{m}_{j=i}\var(Y_{j}-Y_{j-1}).
\end{displaymath}Let \(\tilde\tau^{(m)}\) be the expected time to simulate \(\tilde Z^{(m)}\).  Then\begin{equation*}
E(\tilde\tau^{(m)})=\sum^{m}_{i=0}q_{i}(\tilde t_{i+1}-\tilde t_{i}).
\end{equation*}
Thus,  \(q^{*}=\arg\min_{q\in A_{m}}R^{(m)}(q;\tilde t^{(m)},\tilde\eta)\) optimizes the performance of   \( \tilde Z^{(m)}\).  The vector \(q^{*}\) can be found in \(O(m)\) time by replacing \((\bar t^{(m)},\bar\eta)\) with \((\tilde t^{(m)},\tilde\eta)\)   in the algorithm described in \S\ref{sub:hull}. Here again, the vector   \(\tilde\eta\)  can be estimated by Monte Carlo simulation. 
\section{Numerical Experiments}\label{se:numer}
The simulation experiments were implemented  in the C++ programming language.
The price \(X(t)\) of a stock at time \(t\) is assumed to follow the Geometric Brownian motion model (e.g. \citet{glasserman2004Monte}), i.e. it satisfies the stochastic differential equation
\begin{equation*}
dX(t)= rX(t)\,dt+\sigma X(t)\,dW,
\end{equation*} with \(X(0)=1\), where \(r=0.05\) is the risk-free rate, \(\sigma=0.2\) is the volatility of the stock, and \(W\) is a Brownian motion under the risk-neutral probability. Here \(Y=f(X(T))\), where \(T=1\) and \(f(x)=e^{-r}\max(x-1,0)\), so that \(E(Y)\) is the price of a one-year at the money call option on the stock. The model parameters are the same as those in~\cite{GlynnRhee2015unbiased}, and \(E(Y)\) is approximately equal to \(0.104505836\). The numerical experiments approximate \(Y\) via the Milstein scheme described in Example~\ref{ex:option}. For the coupled sum (resp. independent sum) estimator, \(q_{0},\dots,q_{m}\)  are calculated via the algorithm described in \S\ref{sub:hull} (resp. \S\ref{sub:trIndSum}), with \(m=13\). Each component of the vectors \(\bar\eta\) and \(\tilde\eta\) is estimated via Monte Carlo simulation using \(10000\) independent runs. The remaining components  of \(q\) are calculated by setting \(q_{i}=2^{-3(i-m)/2}q_{m}\), \(i\ge m+1\).  In Table~\ref{tab:BS}, the estimated option price and its standard deviation Std  are calculated using \(n\) independent copies of the corresponding estimator.
 The variable Work refers to the total expected number of time steps  simulated through the \(n\) replications. The
fifth column is a \(90\%\)-confidence interval for this variable.
 Thus, \(\text{Work} \times \text{Std}^2$ is an estimate of the work-normalized variance. As expected, in Table~\ref{tab:BS},  for each of the coupled  and independent sum estimators, the variable Work is roughly proportional to \(n\), the variable Std is roughly proportional to \(n^{-1/2}\), while the work-normalized variance is roughly independent of \(n\). The two estimators have a similar work-normalized variance, and thus a similar performance.\begin{table}
\caption{Pricing of a call option under a Geometric Brownian motion model.}
\begin{small}
\begin{tabular}{lrlccc}\hline
    Estimator & $n$& price & Std & Work& Work $\times$ Std$^2$  \\ \hline
Coupled Sum & $ 10^4$ & $0.1063 $ & $1.7\times10^{-3}$ & $1.108\times10^4 \pm 1.7\times10^2$ & $0.030$\\
Coupled Sum & $ 10^6$ & $0.10445 $ & $1.6\times10^{-4}$ & $1.124\times10^6 \pm 4.9\times10^3$ & $0.029$\\
Coupled Sum & $ 10^{8}$ & $0.104505 $ & $1.6\times10^{-5}$ & $1.125\times10^8 \pm 9.8\times10^4$ & $0.031$\\
Coupled Sum & $ 10^{10}$ & $0.1045066 $ & $1.7\times10^{-6}$ & $1.126\times10^{10} \pm 9.5\times10^5$ & $0.031$\\
Independent Sum & $ 10^4$ & $0.1042 $ & $1.6\times10^{-3}$ & $1.201\times10^4 \pm 1.7\times10^2$ & $0.029$\\
Independent Sum & $ 10^6$ & $0.10458 $ & $1.5\times10^{-4}$ & $1.191\times10^6 \pm 4.9\times10^3$ & $0.028$\\
Independent Sum & $ 10^{8}$ & $0.104521 $ & $1.5\times10^{-5}$ & $1.191\times10^8 \pm 9.8\times10^4$ & $0.028$\\
Independent Sum & $ 10^{10}$ & $0.1045064 $ & $1.5\times10^{-6}$ & $1.192\times10^{10} \pm 9.5\times10^5$ & $0.031$\\
\hline\end{tabular}
 \end{small}
 \label{tab:BS}
\end{table}

\section{Conclusion}\label{se:conc}This note establishes necessary and sufficient conditions for the square integrability of the coupled sum and independent sum estimators.  These conditions are weaker than the sufficient condition of \citet{GlynnRhee2015unbiased}. A geometric characterization of a distribution with infinite support that optimizes the performance of these estimators is presented. An algorithm based on convex hulls that finds an optimal \(m\)-truncated distribution in \(O(m)\) time is described. The algorithm is simple to implement and is illustrated using a numerical example. Alternative RMLMC estimators not covered in this note, such as the ``single term'' estimator, are studied in \cite{GlynnRhee2015unbiased,Vihola2018}. Using the results in this note to  broaden the range of applications of the RMLMC methods is a promising direction for future research. 

\commentt{\appendix}
{
\begin{APPENDIX}{}
}
\section{Proof of Theorem~\ref{th:GenCase}}
Let us first prove the following.
\begin{proposition}\label{pr:nuq}
Let \((\gamma_{n})\), \(n\ge0\), be a nonnegative sequence that goes to \(0\) as \(n\) goes to infinity. Then there is a strictly increasing nonnegative  integral sequence  \((\rho(n):n\ge0)\),  such that
\begin{equation}\label{eq:rhoCond}  
\gamma_{\rho(n)}\le\gamma_{j}\text{ for }n\ge0\text{ and }j\in[\rho(0),\rho(n)].
\end{equation}      
Furthermore, if \(q\in A\) and\begin{equation}\label{eq:NewCondZIntegralNu}
\sum^{\infty}_{i=0}(\frac{1}{q_{i+1}}-\frac{1}{q_{i}})\gamma_{i}<\infty,
\end{equation} then \(\gamma_{\rho(n)}/q_{\rho(n)+1}\) goes to \(0\) as \(n\) goes to infinity. 
\end{proposition}
\commentt{\begin{proof}}{\proof{Proof.}}
Assume first that there  is an integer \(n_{0}\) such that  \(\gamma_{n}>0\) for \(n\geq n_{0}\).  Let \(\rho(0)=n_0\) and, for \(n\ge1\), let\begin{displaymath}
\rho(n)=\min\{j > \rho(n-1):\gamma_{j}\le \gamma_{\rho(n-1)}\}.
\end{displaymath}By construction, \eqref{eq:rhoCond} holds. Assume now that there are infinitely many integers \(n\) such that \(\gamma_{n}=0\). Let \(\rho\) be a strictly increasing sequence with \(\gamma_{\rho(n)}=0\) for \(n\ge0\). Here again,  \eqref{eq:rhoCond} holds. 
 
Suppose now that  \(q\in A\) and that \eqref{eq:NewCondZIntegralNu} holds. Fix \(\epsilon>0\). By~\eqref{eq:NewCondZIntegralNu}, there is an integer \(m\) such that, for \(n>m\), \begin{equation*}
\sum^{\rho(n)}_{i=\rho(m)}(\frac{1}{q_{i+1}}-\frac{1}{q_{i}})\gamma_{i}<\epsilon/2,
\end{equation*}
and so, by  \eqref{eq:rhoCond}, \begin{displaymath}
(\frac{1}{q_{\rho(n)+1}}-\frac{1}{q_{\rho(m)}})\gamma_{\rho(n)}<\epsilon/2.
\end{displaymath} 
Because the sequence  \(q\) goes to \(0\) at infinity, there is an integer \(m'>m\)   such that   \(q_{\rho(n)+1}<q_{\rho(m)}/2\) for \(n>m'\). Thus  \(\gamma_{\rho(n)}/q_{\rho(n)+1}<\epsilon\) for \(n>m'\). This concludes the proof.\commentt{\end{proof}}{\Halmos\endproof}
\comment{\begin{lemma}Let \((\Delta_{n})\), \(n\geq0\), be a sequence of square-integrable random variables. For \(n\geq0\),  set
 \begin{displaymath}
V_{n}=\sum^{n\wedge N}_{k=0}\frac{\Delta _{k}}{q_{k}}.
\end{displaymath}
Let \(V=\sum^{N}_{k=0}{\Delta _{k}}/{q_{k}}\).   
\end{lemma}
}
For the rest of the paper, for \(n\ge0\), set \(\Delta_{n}=Y_{n}-Y_{n-1}\) and
 \begin{displaymath}
\bar Z_{n}=\sum^{n\wedge N}_{k=0}\frac{\Delta _{k}}{q_{k}},
\end{displaymath}
with \(\bar Z_{-1}=0\). Proposition~\ref{pr:RheeGlynnBarZ} is shown in \cite[pp. 1027, 1030]{GlynnRhee2015unbiased}.  
\begin{proposition}\label{pr:RheeGlynnBarZ}[\citet{GlynnRhee2015unbiased}]
For \(-1\le m\le n\),
we have \(E(\bar Z_{n})=E(Y_{n})\), and \begin{equation*}
||\bar Z_{n}-\bar Z_{m}||^{2}=\sum^{ n}_{i= m+1}\frac{||Y_{i-1}-Y_{ n}||^{2}-||Y_{i}-Y_{ n}||^{2}}{q_{i}}
.
\end{equation*}\end{proposition}
 
\begin{proposition}
\label{pr:subSequenceBounded}If \(\bar Z\) is square-integrable then    \(||\bar Z_{ n}{\mathbf 1}\{N>  n\}||^{2}\leq||\bar Z||^{2}+1\) for infinitely many  integers  \(n\). 
\end{proposition}
\commentt{\begin{proof}}{\proof{Proof.}}
Fix \(n\ge0\) with \(\Pr(N=n)>0\). Because of the independence of \(N\) and \((Y_{i}:i\ge0)\),\begin{eqnarray*}E({\bar Z_{n} }^2|N>n)&=&||\sum^{n}_{k=0}\frac{\Delta _{k}}{q_{k}}||^{2}\\
&=& E({\bar Z}^2|N=n).
\end{eqnarray*}Hence 
\begin{equation}\label{eq:prbarZn}
||\bar Z_{n}{\mathbf 1}\{N>n\}||^{2}\Pr(N=n)=||\bar Z{\mathbf 1}\{N=n\}||^{2}
\Pr(N>n).
\end{equation}
Clearly, this equation also holds if  \(\Pr(N=n)=0\). 

Assume now for contradiction that the conclusion of Proposition~\ref{pr:subSequenceBounded} does not hold. Then there is an integer \(m\)  such that  \(||\bar Z_{ n}{\mathbf 1}\{N>  n\}||^{2}>||\bar Z||^{2}+1\) for \(n\ge m\). By \eqref{eq:prbarZn}, for \(n\ge m\),

\begin{equation*}
||\bar Z_{n}{\mathbf 1}\{N>n\}||^{2}\Pr(N=n)\le q_{m}||\bar Z{\mathbf 1}\{N=n\}||^{2},
\end{equation*}
and so  
 
\begin{equation*}
(||\bar Z||^{2}+1)\Pr(N=n)\le q_{m}||\bar Z{\mathbf 1}\{N=n\}||^{2}.
\end{equation*}Summing over \(n\in[m,\infty)\) implies that

\begin{equation*}
(||\bar Z||^{2}+1)q_{m}\le q_{m}||\bar Z||^{2},
\end{equation*} leading to a contradiction. \commentt{\end{proof}}{\Halmos\endproof}
Let us now prove Theorem~\ref{th:GenCase}. The first part of the proof is inspired from the proof of Theorem~1 of \cite{GlynnRhee2015unbiased}. By Proposition~\ref{pr:RheeGlynnBarZ},  for \(-1\le m\le n\),
\begin{equation}\label{eq:thBarZrewritten}
||\bar Z_{n}-\bar Z_{m}||^{2}=\frac{||Y_{m}-Y_{n}||^{2}}{q_{m+1}}
+\sum^{n}_{i=m+1}||Y_{i}-Y_{n}||^{2}(\frac{1}{q_{i+1}}-\frac{1}{q_{i}}). \end{equation} 

Suppose first that~\eqref{eq:NewCondZIntegral} holds. Applying Proposition~\ref{pr:nuq} with \(\gamma_{n}=||Y_{n}-Y||^{2}\) shows the existence of a  strictly increasing nonnegative  integral sequence  \((\rho(n):n\ge0)\) such that  \eqref{eq:rhoCond} holds.         Set \(\bar Z'_{n}=\bar Z_{\rho(n)}\) for \(n\geq0\).  By \eqref{eq:thBarZrewritten},  for \(0\le m\le n\),
\begin{equation*}
||\bar Z'_{n}-\bar Z'_{m}||^{2}=\frac{||Y_{\rho(m)}-Y_{\rho(n)}||^{2}}{q_{\rho(m)+1}}
+\sum^{\rho(n)}_{i=\rho(m)+1}||Y_{i}-Y_{\rho(n)}||^{2}(\frac{1}{q_{i+1}}-\frac{1}{q_{i}}). \end{equation*}
Let us show that \((\bar Z'_{n}:n\ge0)\) is a Cauchy sequence in \(L^{2}\). For \(\rho({0})\le i\le\rho(n)\),
\begin{eqnarray}\label{eq:upperBoundYi}
||Y_{i}-Y_{\rho(n)}||^{2}&\le&2(\gamma_{i}+\gamma_{\rho(n)})\nonumber\\
&\le&4\gamma_{i},
\end{eqnarray}
where the second inequality follows from   \eqref{eq:rhoCond}. Thus, for \(0\leq m\leq n\),
\begin{equation}\label{eq:Cauchy}
||\bar Z'_{n}-\bar Z'_{m}||^{2}\le 4\frac{\gamma_{\rho(m)}}{q_{\rho(m)+1}}
+4\sum^{\rho(n)}_{i=\rho(m)+1}\gamma_{i}(\frac{1}{q_{i+1}}-\frac{1}{q_{i}}). \end{equation} By   \eqref{eq:NewCondZIntegral} and Proposition~\ref{pr:nuq},  for any \(\epsilon>0\), the first term in the RHS of~\eqref{eq:Cauchy}
is smaller than \(\epsilon\) if \(m\) is sufficiently large. Because of~\eqref{eq:NewCondZIntegral}, the same holds for the second term. Thus, the sequence \((\bar Z'_{n}:n\ge0)\) is Cauchy in \(L^{2}\), and so it has a limit in \(L^{2}\) as \(n\) goes to infinity. Since \(\bar Z'_{n}\) converges a.s. to \(\bar Z\) as \(n\) goes to infinity, this implies  that  \(\bar Z\) is in \(L^{2}\) and that \(\bar Z'_{n}\) converges in \(L^{2}\) to \(\bar Z\) as \(n\) goes to infinity. Hence \(E(\bar Z'_{n})\) (resp.   \(||\bar Z'_{n}||\)) converges to \(E(\bar Z)\) (resp. \(||\bar Z||\)) as \(n\) goes to infinity. By Proposition \ref{pr:RheeGlynnBarZ}, \(E(\bar Z'_{n})=E(Y_{\rho(n)})\). Letting \(n\) go to infinity implies that \(E(\bar Z)=E(Y)\). This is because \(Y_{n}\) converges to \(Y\) in \(L^{2}\). Furthermore, applying \eqref{eq:thBarZrewritten} with \(m=-1\) and replacing \(n\) by \(\rho(n)\) yields 
\begin{equation*}
||\bar Z'_{n}||^{2}=||Y_{\rho(n)}||^{2}+\sum^{\rho(n)-1}_{i=0}||Y_{i}-Y_{\rho(n)}||^{2}(\frac{1}{q_{i+1}}-\frac{1}{q_{i}}).
\end{equation*}
 Observe that  \(||Y_{\rho(n)}||\) (resp. \(||Y_{i}-Y_{\rho(n)}||\)) converges to \(||Y||\) (resp. \(||Y_{i}-Y||\)) as \(n\) goes to infinity. Letting \(n\) go to infinity and using \eqref{eq:upperBoundYi}, \eqref{eq:NewCondZIntegral}, and 
    the dominated convergence theorem yields~\eqref{eq:normbarZ}.

Assume now  that \(\bar Z\) is square-integrable. For \(n\ge0\),  by linearity of expectation and the equality   \(\bar Z_{n}{\mathbf 1}\{N\leq n\}=\bar Z{\mathbf 1}\{N\leq n\}\), 
\begin{eqnarray*}
||\bar Z_{n}||^{2}&=&||\bar Z{\mathbf 1}\{N\leq n\}||^{2}+||\bar Z_{n}{\mathbf 1}\{N> n\}||^{2}\\&\le&||\bar Z||^{2}+||\bar Z_{n}{\mathbf 1}\{N> n\}||^{2}.
\end{eqnarray*}
Combining this with Proposition~\ref{pr:subSequenceBounded} shows the existence of a strictly increasing nonnegative integral sequence  \(\lambda\) such that, for \(n\in\mathbb{N}\),
\begin{equation*}
||\bar Z_{\lambda(n)}||^{2}\le2||\bar Z||^{2}+1.
\end{equation*}   
Applying once again  \eqref{eq:thBarZrewritten} with \(m=-1\) yields
\begin{equation*}
||Y_{\lambda(n)}||^{2}+\sum^{\lambda(n)-1}_{i=0}||Y_{i}-Y_{\lambda(n)}||^{2}(\frac{1}{q_{i+1}}-\frac{1}{q_{i}})\le2||\bar Z||^{2}+1.
\end{equation*}
For \(m\in\mathbb{N}\) and  \(n\ge m\), because \((q_{i})\) is a decreasing sequence and \(m\le\lambda(n)\), it follows that
\begin{equation*}
\sum^{ m-1}_{i=0}||Y_{i}-Y_{\lambda(n)}||^{2}(\frac{1}{q_{i+1}}-\frac{1}{q_{i}})\le2||\bar Z||^{2}+1.
\end{equation*}
Letting \(n\) go to infinity shows that \begin{equation*}
\sum^{ m-1}_{i=0}||Y_{i}-Y||^{2}(\frac{1}{q_{i+1}}-\frac{1}{q_{i}})\le2||\bar Z||^{2}+1,
\end{equation*}
which implies~\eqref{eq:NewCondZIntegral}.\qed   
\section{Proof of Theorem~\ref{th:GenCaseTildeZ}}
The proof is similar to that of  Theorem~\ref{th:GenCase}. For \(n\geq0\), set
 \begin{displaymath}
\tilde Z_{n}=\sum^{n\wedge N}_{k=0}\frac{\tilde\Delta _{k}}{q_{k}},
\end{displaymath}
and let \(\tilde Z_{-1}=0\). Note that \(E(\tilde Z_{n})=E(\bar Z_{n})=E(Y_{n})\). Proposition~\ref{pr:RheeGlynnTildeZ} is shown in  \cite[p. 1031]{GlynnRhee2015unbiased}.  
\begin{proposition}\label{pr:RheeGlynnTildeZ}[\citet{GlynnRhee2015unbiased}]
For \(-1\le m\le n\),
\begin{equation*}
||\tilde Z_{n}-\tilde Z_{m}||^{2}=\sum^{ n}_{i= m+1}\frac{E(\Delta_i^{2}+2\Delta_{i}E(Y_{n}-Y_{ i}))}{q_{i}}
.
\end{equation*}\end{proposition}
  A simple calculation shows that \begin{displaymath}
E(\Delta_i^{2}+2\Delta_{i}E(Y_{n}-Y_{ i}))=\var(\Delta_{i})+(E(Y_{n}-Y_{ i-1}))^{2}-(E(Y_{n}-Y_{ i}))^{2}.
\end{displaymath} Hence, by Proposition~\ref{pr:RheeGlynnTildeZ},  for \(-1\le m\le n\),
\begin{equation}\label{eq:thTildeZrewritten}
||\tilde Z_{n}-\tilde Z_{m}||^{2}=\frac{(E(Y_{m}-Y_{n}))^{2}}{q_{m+1}}
+\sum^{ n}_{i= m+1}\left(\frac{\var(\Delta_i)}{q_{i}}
+(E(Y_{i}-Y_{n}))^{2}(\frac{1}{q_{i+1}}-\frac{1}{q_{i}})\right). \end{equation} 

Suppose first that~\eqref{eq:NewCondZIntegralTilde} holds. Applying Proposition~\ref{pr:nuq} with \(\gamma_{n}=(E(Y_{n}-Y))^{2}\) shows the existence of a strictly increasing nonnegative   integral sequence  \((\rho(n):n\ge0)\) such that  \(|E(Y_{\rho(n)}-Y)|\le|E(Y_{j}-Y)|\) for \(\rho(0)\leq j\leq \rho(n)\).         Set \(\tilde Z'_{n}=\tilde Z_{\rho(n)}\) for \(n\geq0\).  By \eqref{eq:thTildeZrewritten},  for \(0\le m\le n\),
 \begin{equation*}
||\tilde Z'_{n}-\tilde Z'_{m}||^{2}=\frac{(E(Y_{\rho(m)}-Y_{\rho(n)}))^{2}}{q_{\rho(m)+1}}
+\sum^{\rho(n)}_{i= \rho(m)+1}\left(\frac{\var(\Delta_i)}{q_{i}}
+(E(Y_{i}-Y_{\rho(n)}))^{2}(\frac{1}{q_{i+1}}-\frac{1}{q_{i}})\right). 
\end{equation*}
For \(\rho({0})\le i\le\rho(n)\),
\begin{eqnarray*}
|E(Y_{i}-Y_{\rho(n)})|&\le&|E(Y_{i}-Y)|+|E(Y_{\rho(n)}-Y)|\nonumber\\
&\le&2|E(Y_{i}-Y)|.
\end{eqnarray*}
Thus, for \(0\leq m\leq n\),
\begin{equation}\label{eq:CauchyTilde}
||\tilde Z'_{n}-\tilde Z'_{m}||^{2}\le 4\frac{\gamma_{\rho(m)}}{q_{\rho(m)+1}}
+\sum^{\rho(n)}_{i=\rho(m)+1}\left(\frac{\var(\Delta_i)}{q_{i}}+4\gamma_{i}(\frac{1}{q_{i+1}}-\frac{1}{q_{i}})\right). \end{equation} By   \eqref{eq:NewCondZIntegralTilde} and Proposition~\ref{pr:nuq},  for any \(\epsilon>0\), the first term in the RHS of~\eqref{eq:CauchyTilde}
is smaller than \(\epsilon\) if \(m\) is sufficiently large. Because of~\eqref{eq:NewCondZIntegralTilde}, the same holds for the second term. Thus, the sequence \((\tilde Z'_{n}:n\ge0)\) is  Cauchy  in \(L^{2}\).   As in the proof   Theorem~\ref{th:GenCase}, this implies that  \(\tilde Z\) is in \(L^{2}\), that  \(E(\tilde Z)=E(Y)\), and that   \(||\tilde Z'_{n}||\) converges to  \(||\tilde Z||\) as \(n\) goes to infinity.  Furthermore, applying \eqref{eq:thTildeZrewritten} with \(m=-1\) and replacing \(n\) by \(\rho(n)\) yields
\begin{equation*}
||\tilde Z'_{n}||^{2}=(E(Y_{\rho(n)}))^{2}
+\sum^{ \rho(n)}_{i=0}\left(\frac{\var(\Delta_i)}{q_{i}}
+(E(Y_{i}-Y_{\rho(n)}))^{2}(\frac{1}{q_{i+1}}-\frac{1}{q_{i}})\right). 
\end{equation*}
  Letting \(n\) go to infinity and using the dominated convergence theorem implies~\eqref{eq:normTildeZ}.

Conversely, if \(\tilde Z\) is square-integrable, then    \eqref{eq:NewCondZIntegralTilde} follows from arguments similar to those used in the proof of  Theorem~\ref{th:GenCase}.\qed
\section{Proof of Proposition \ref{pr:varZtilde}}
Let us first prove the following proposition. The equality \eqref{eq:intPartInfinite}  means that, if one of its members is finite, so is the other one, and the two members are equal. Moreover, if one of its members is infinite, so  is the other one.  
\begin{proposition}\label{pr:nuqFirst}
Let \((\gamma_{n}:n\ge0)\) be a positive decreasing sequence that goes to \(0\) as \(n\) goes to infinity, and let \(q\in A\). Then \begin{equation}\label{eq:intPartInfinite}
\sum^{\infty}_{i=0}\frac{\gamma_{i}-\gamma_{i+1}}{ q_{i}}=\gamma_{0}+\sum^{\infty}_{i=1}\gamma _{i}(\frac{1}{q_{i}}-\frac{1}{q_{i-1}}).
\end{equation}
\end{proposition}
\commentt{\begin{proof}}{\proof{Proof.}}
For \(n\geq0\), \begin{equation}\label{eq:intPartFinite}
\sum^{n}_{i=0}\frac{\gamma_{i}-\gamma_{i+1}}{ q_{i}}=\gamma_{0}-\frac{\gamma _{n+1}}{q_{n+1}}+\sum^{n+1}_{i=1}\gamma _{i}(\frac{1}{q_{i}}-\frac{1}{q_{i-1}}).
\end{equation}
Assume first that the LHS of \eqref{eq:intPartInfinite} is infinite. Then  \eqref{eq:intPartFinite} implies that the RHS of \eqref{eq:intPartInfinite} is infinite as well. Assume now that the LHS of \eqref{eq:intPartInfinite} is finite. As \(q\) is decreasing, for \(0\leq m< n\),
\begin{equation*}
\sum^{n}_{i=m}\frac{\gamma_{i}-\gamma_{i+1}}{q_{i}}\ge\frac{\gamma_{m}-\gamma_{n+1}}{ q_{m}}.
\end{equation*}
Letting \(n\) go to infinity shows that 
\begin{equation*}
\frac{\gamma_{m}}{ q_{m}}\le\sum^{\infty}_{i=m}\frac{\gamma_{i}-\gamma_{i+1}}{q_{i}},
\end{equation*} and so \({\gamma_{m}}/{ q_{m}}\) goes to \(0\) as \(m\) goes to infinity. Letting \(n\) go to infinity in \eqref{eq:intPartFinite} implies \eqref{eq:intPartInfinite}.
\commentt{\end{proof}}{\Halmos\endproof}
Let us now prove Proposition \ref{pr:varZtilde}. Assume that    \(q\in A\) is such that \(\tilde Z\) is square-integrable.  Since \(q_{i}\leq1\) for \(i\geq0\), by Theorem~\ref{th:GenCaseTildeZ},  \(\tilde\mu_{0}\) is finite. For \(i\geq0\), let \begin{displaymath}
\gamma_{i}=\sum^{\infty}_{j=i}\var(Y_{j}-Y_{j-1}).
\end{displaymath}
 By  Theorem~\ref{th:GenCaseTildeZ}  and Proposition \ref{pr:nuqFirst},  
\begin{eqnarray*}
\var(\tilde Z)
&=&\sum^{ \infty}_{i=0}\left(\frac{\gamma_{i}-\gamma_{i+1}}{q_{i}}
+(\frac{1}{q_{i+1}}-\frac{1}{q_{i}})(E(Y_{i}-Y))^{2}\right)\\
&=&\gamma_{0}+\sum^{\infty}_{i=1}(\frac{1}{q_{i}}-\frac{1}{q_{i-1}})\gamma_{i}+\sum^{ \infty}_{i=1}(\frac{1}{q_{i}}-\frac{1}{q_{i-1}})(E(Y_{i-1}-Y))^{2}.\end{eqnarray*}
This concludes the proof.\qed
\section{Proof of Theorem~\ref{th:OptTimeVarianceGen}}

The proof builds on ideas used in  the proof of \cite[Theorem 3]{kahaRandomizedDimensionReduction19}. It  uses the  following proposition, whose proof follows immediately from \eqref{eq:RDefGen}.

\begin{proposition}
\label{pr:lowerW}Let  \(\vartheta=(\vartheta_{0},\vartheta_{1},\dots)\in\{0\}\times\mathbb{R}_{+}^{\mathbb{N}}\) be a  strictly increasing sequence. If   \(\gamma\in\mathbb{R}_{+}^{\mathbb{N}}\) and \(\gamma'\in\mathbb{R}_{+}^{\mathbb{N}}\)
are such that \(\gamma'\leq\gamma\),  and  \(q\in A\). If  \(R(q;\vartheta,\gamma)\)is finite then \(R(q;\vartheta,\gamma')\leq R(q;\vartheta,\gamma)\),
 with equality if  \(\gamma_{0}=\gamma'_{0}\) and, for \(i\geq1\),  \((\gamma_{i}-\gamma'_{i})(q_{i-1}- q_{i})=0\). \end{proposition}
 Let us show that \(q^{*}\) is well-defined and belongs to \(A\).
 \comment{ Define the sequence \((i_{k})\), \(k\geq0\), recursively as follows. Set \(i_{0}=0\). Let \(k\geq0\). As \(\gamma_{j}\) goes to \(0\) as \(j\) goes to infinity, the sequence \((\gamma_{j}-\gamma_{i_{k}})/(\bar t_{j}-\bar t_{i_{k}})\), \(j>i_{k}\), contains a negative element and goes to \(0\) as \(j\) goes to infinity. Let \(i_{k+1}\) be the smallest index \(j\) where this sequence attains its minimum. Thus \begin{equation}
i_{k+1}= \arg\min_{j>i_{k}}\frac{\gamma_{j}-\gamma_{i_{k}}}{\bar t_{j}-\bar t_{i_{k}}}.
\end{equation}}Fix \(i\geq 0\), and let  \begin{displaymath}
a_{i} =\min_{0\leq k\leq i}(\frac{\gamma _{k}}{\vartheta_{i+1}-\vartheta_{k}}).
\end{displaymath}  The sequence \((a_{i}(\vartheta_{i+1}-\vartheta_{n}):n\geq0)\) is an affine function of  \(\vartheta\) and is upper-bounded by \(\gamma\), and so it is upper-bounded by \(\gamma'\). Hence \(\gamma'_i\ge a_{i}(\vartheta_{i+1}-\vartheta_{i})>0\). In the particular case where \(i=0\), this implies that \(\gamma'_{0}=\gamma_{0}\).  Since    \(\gamma_{n}\)
goes to \(0\) as \(n\) goes to infinity, there is \(j>i \) with \(\gamma_j<\gamma'_{i}\). By definition of the lower hull,  \(\gamma_j\ge \gamma'_j\ge\gamma'_i+\theta_{i}(\vartheta_{j}-\vartheta_{i})\), and so
\(\theta_{i}< 0\). Thus   \(q^{*}_{i}\) is well-defined and is   strictly positive.
On the other hand, for \(n>0\), because of the convexity properties of the lower hull, \begin{displaymath}
\theta_{n}\geq\frac{\gamma'_{n}-\gamma'_{0}}{\vartheta_{n}}\ge\frac{-\gamma_{0}}{\vartheta_{n}} .
\end{displaymath} Thus \(|\theta_{n}|\le\gamma_{0}/\vartheta_{n}\), and so \(\theta_{n}\) goes to \(0\) as \(n\) goes to infinity. As the sequence \((\theta_{n}:n\geq0)\) is increasing, this implies that  \(q^{*}\in A\). Furthermore, since 
\(\theta_{n}< 0\) for \(n\ge0\), the  sequence   \(\gamma'\) is strictly decreasing.

By hypothesis, there is \(q\in A\) such that \(R(q;\vartheta,\gamma)\) is finite. Thus, each  component of the product in the RHS of \eqref{eq:RDefGen} has a finite limit as \(n\) goes to infinity.  By Proposition~\ref{pr:lowerW},   \(R(q;\vartheta,\gamma')\le R(q;\vartheta,\gamma)\), and so  \(R(q;\vartheta,\gamma')\) is finite as well. Since \(\gamma'\) is positive, decreasing and goes to \(0\) at infinity, Proposition~\ref{pr:nuqFirst} shows that
\begin{equation*}
\sum^{\infty}_{i=0}\frac{\gamma'_{i}-\gamma'_{i+1}}{ q_{i}}=\gamma'_{0}+\sum^{\infty }_{i=1}\gamma' _{i}(\frac{1}{q_{i}}-\frac{1}{q_{i-1}}).
\end{equation*} As \(\gamma'\le\gamma\), the RHS of this equation is finite. Thus,  
   \begin{equation*}
R(q;\vartheta,\gamma')=(\sum^{\infty}_{i=0}q_{i}(\vartheta_{i+1}-\vartheta_{i}))(\sum^{\infty}_{i=0}\frac{\gamma'_{i}-\gamma'_{i+1}}{ q_{i}}),
\end{equation*} and the  two series in the RHS of this equation are convergent. Since, by the Cauchy-Schwartz inequality, for all nonnegative summable sequences \((x_{i})\) and \((y_{i})\), \(i\ge0\),\begin{displaymath}
(\sum^{\infty}_{i=0}\sqrt{x_{i}y_{i}})^{2}\le(\sum^{\infty}_{i=0}x_{i})(\sum^{\infty}_{i=0}y_{i}),
\end{displaymath}
we conclude that\begin{equation}\label{eq:barRtheoProof}
\bigg(\sum_{i=0}^{\infty}\sqrt{(\gamma'_{i}-\gamma'_{i+1})(\vartheta_{i+1}-\vartheta_{i})}\bigg
)^2\leq R(q;\vartheta,\gamma'),
\end{equation} 
and that the LHS of this equation is finite. Furthermore, by the definition of \(q^{*}\) and Proposition~\ref{pr:nuqFirst}, \begin{eqnarray*}\sqrt{|\theta_{0}|}\sum_{i=0}^{\infty}\sqrt{(\gamma'_{i}-\gamma'_{i+1})(\vartheta_{i+1}-\vartheta_{i})}&=& \sum^{\infty}_{i=0}\frac{\gamma'_{i}-\gamma'_{i+1}}{ q^{*}_{i}}\\
\\&=&\gamma'_{0}+\sum^{\infty}_{i=1}(\frac{1}{q^{*}_{i}}-\frac{1}{q^{*}_{i-1}})\gamma'_{i}.\end{eqnarray*}
The second equation follows once again from  Proposition~\ref{pr:nuqFirst}. Thus the three series are convergent and have the same limit. Using~\eqref{eq:RDefGen} and the definition of \(q^{*}\) implies that
\begin{equation*}
R(q^{*};\vartheta,\gamma')=\bigg(\sum_{i=0}^{\infty}\sqrt{(\gamma'_{i}-\gamma'_{i+1})(\vartheta_{i+1}-\vartheta_{i})}\bigg
)^2.
\end{equation*} Hence, by~\eqref{eq:barRtheoProof}, \(R(q^{*};\vartheta,\gamma')\leq R(q;\vartheta,\gamma')\).  On the other hand,  \((\gamma_{i}-\gamma'_{i})(q^{*}_{i-1}- q^{*}_{i})=0\) for \(i\geq 1\). This is because, if  \(\gamma'_{i}<\gamma_{i}\), then the point \((\vartheta_{i},\gamma_{i})\) does not belong to the lower hull of the set  \(\{(\vartheta_{j},\gamma_{j}):j\geq 0\}\). Thus  \((\vartheta_{i},\gamma'_{i})\) belongs to the segment \([ (\vartheta_{i-1},\gamma'_{i-1}),(\vartheta_{i+1},\gamma'_{i+1})] \), which implies that \(\theta_{i-1}= \theta_{i}\) and \(q^{*}_{i-1}= q^{*}_{i}\).   Hence, as   \(\gamma_{0}=\gamma'_{0}\),   Proposition~\ref{pr:lowerW} shows that    \(R(q^{*};\vartheta,\gamma)=\ R(q^{*};\vartheta,\gamma')\).  This implies \eqref{eq:optTimeVariance} and shows that  \(R(q^{*};\vartheta,\gamma)\leq R(q;\vartheta,\gamma)\), as desired.
\commentt{\qed}{\Halmos}
\section*{Acknowledgments}
This work was achieved through the Laboratory of Excellence on
Financial Regulation (Labex ReFi) under the reference ANR-10-LABX-0095. It benefitted from
a French government support managed by the National Research Agency (ANR).

\commentt{}
{
\end{APPENDIX}
}{}
\bibliography{poly}
\end{document}